\hfuzz=9.5pt
\hsize=12.5cm
\vsize=19.2cm
\parindent=0.5cm
\parskip=0pt
\baselineskip=12pt
\topskip=12pt
\magnification=\magstep1

%				NEWLINE
\def\qed{\hbox{\hskip 6pt\vrule width6pt height7pt depth1pt \hskip1pt}}
\def\natural{{\rm I\kern-.18em N}}

\def\integer{{\rm Z\kern-.32em Z}}
\def\chix{{\raise.5ex\hbox{$\chi$}}}

\def\real{{\rm I\kern-.2em R}}

\def\complex{\kern.1em{\raise.47ex\hbox{
	    $\scriptscriptstyle |$}}\kern-.40em{\rm C}}

\def\undertext#1{$\underline{\smash{\hbox{#1}}}$}

%%Margaret's macros for making boxes
%
\def\nl{\par}
\long\def\boxit#1#2{\vbox{\hrule\hbox{\vrule\kern#1
        \vbox{\kern#1\vbox{#2}\kern#1}\kern#1\vrule}\hrule}}

\def\vs#1 {\vskip#1truein}
\def\hs#1 {\hskip#1truein}
  \hbadness=10000 \vbadness=10000 %	REPORT ONLY BEYOND THIS BADNESS
%\nopagenumbers
\def\cala{{\cal A}}
\def\a{\alpha}

\def\nd{\noindent}
\def\calb{{\cal B}}
\def\calg{{\cal G}}
\def\calm{{\cal M}}
\def\tcalm{{\tilde \calm}_{k+1}}
\def\l{{\lambda}}

\def\std{{\rm std}}

\pageno=1

\footline{\ifnum\pageno=0\hss\else\hss\tenrm\folio\hss\fi}
\hbox{}
\vskip 1truein\centerline{{\bf SIMPLE TYPE IS NOT A BOUNDARY PHENOMENON}}
\vs.2
\centerline{by}
\vs.2 
\centerline{Lorenzo Sadun${}^1$}
\footnote{}{1\ Department of Mathematics, University of Texas, Austin, TX 78712, USA.  Email: sadun@math.utexas.edu.  
Research supported in part by an NSF Mathematical Sciences 
Postdoctoral Fellowship and Texas ARP Grant 003658-037 \hfil}
\vs.5 \centerline{{\bf Abstract}}
\vs.1 \nd
This is an expository article, explaining recent work by D. Groisser
and myself [GS] on the extent to which the boundary region of moduli
space contributes to the ``simple type'' condition of Donaldson
theory.  The presentation is intended to complement [GS], presenting
the essential ideas rather than the analytical details.  It is shown
that the boundary region of moduli space contributes $6/64$ of the
homology required for simple type, regardless of the topology or
geometry of the underlying 4-manifold.  The simple type condition thus
reduces to a statement about the interior of moduli space, namely that
the interior of the $k+1$st ASD moduli space, intersected with two
representatives of (4 times) the point class, be homologous to 58
copies of the $k$-th moduli space.  This is peculiar, since the only
known embeddings of the $k$-th moduli space into the $k+1$st involve
Taubes patching, and the image of such an embedding lies entirely in
the boundary region.

\vs1 
\centerline{December 1996} 

\vfill\eject 

In this paper I discuss some recent work of David Groisser and myself
[GS] on how the ``simple type'' condition of Donaldson theory is
related to the geometry of the moduli spaces of anti-self-dual
connections over a given 4-manifold. But before I begin, I must answer
the obvious nagging question: {\it Haven't the Seiberg-Witten
equations made all of Donaldson theory obsolete?}  I obviously think
not.  While SW theory has eclipsed much of Donaldson
theory, it has actually made other uses of Donaldson invariants {\it
more} practical.  This paper, I hope, will serve as an example of how
to extract insight from Donaldson theory in the post-Seiberg-Witten
era.

\bigskip

\nd \undertext{\bf A Historical Digression}

\medskip

Once upon a time, the Yang-Mills (YM) equations were proposed and were
studied without regard to topological consequences.  
The motivation was from physics: the equations accurately describe 
physics at the subnuclear scale, and the ``standard model'' is an $SU(3)
\times SU(2) \times U(1)$ gauge theory.
In the 1970s and 1980s, people began to write down solutions to the YM
equations.  Of particular interest were anti-self-dual (ASD)
connections, which automatically satisfy the YM equations. People
began to look at ``moduli spaces'' of ASD connections, both over
manifolds of physical interest ($\real^4$, or equivalently $S^4$) and
over more general 4-manifolds.  In addition to aiding our
understanding of the connections themselves, these moduli spaces were of
interest in their own right, especially over complex manifolds, where
ASD connections correspond to holomorphic vector bundles.

Donaldson's brilliant insight was that the topology of moduli spaces
(and in particular how they sit in the larger space of all connections
modulo gauge transformations) could tell a lot about the differential
topology of the underlying 4-manifolds.  Donaldson invariants capture
the essential topological information about the moduli spaces, and
proved very useful for classifying smooth 4-manifolds.  In fact, they
were {\it so} successful that the previous quest, to understand the YM
equations and the ASD moduli spaces for their own sake, was largely
neglected.

By 1990, then, mathematicians primarily studied the YM
equations in order to understand moduli spaces, studied moduli spaces
in order to understand Donaldson invariants, and studied Donaldson
invariants in order to aid in the classification of smooth 4-manifolds.
This is schematically shown in Figure 1.  This was a grand and
difficult project, with a huge number of practitioners making gradual
progress.  

\medskip

$$
\boxit{.5em}{\multiply\hsize2\divide\hsize16 \noindent Yang-Mills\nl
        \noindent Equations} 
\enspace \raise3ex\hbox{$\Rightarrow$}\enspace
\boxit{.5em}{\multiply\hsize2\divide\hsize24 \noindent Moduli\nl
        \noindent Spaces} 
\enspace \raise3ex\hbox{$\Rightarrow$}\enspace
\boxit{.5em}{\multiply\hsize2\divide\hsize16 \noindent Donaldson\nl
        \noindent Invariants} 
\enspace \raise3ex\hbox{$\Rightarrow$}\enspace
\boxit{.5em}{\multiply\hsize2\divide\hsize15 \noindent Classifying\nl
        \noindent 4-manifolds}
$$

\medskip

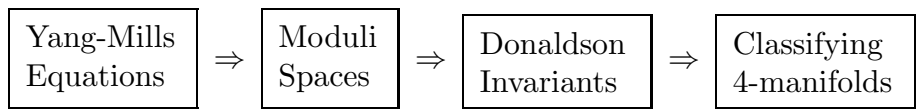
\centerline{Figure 1. The Traditional Flow of Ideas}

%\vfill\eject

\bigskip

\nd \undertext{\bf The Seiberg-Witten Revolution}

\medskip

In late 1994, this project was largely made irrelevant by the advent
of the Seiberg-Witten equations [W].  The SW invariants are far easier to
compute than Donaldson invariants, and are generally believed to carry
exactly the same information.  Witten's formula [W] for the Donaldson
invariants in terms of the SW invariants is almost universally
believed, although as of this writing a mathematical proof is still
lacking.  As far as classifying smooth 4-manifolds goes, just about
anything that can be done with Donaldson theory can be done, far more
easily, with SW invariants.

The last arrow in Figure 1 must therefore be abandoned.  Indeed, the
middle arrow is also largely superceded, as (modulo a proof of the
Witten conjecture) the Donaldson invariants are most easily computed
by first computing the SW invariants and then applying Witten's
formula. The flow of ideas in Figure 1 is simply obsolete.

$SU(2)$ gauge theory, however, is {\it not} obsolete, if we
remember the original interest in the YM equations and the moduli
spaces.  We just have to reverse the arrows in Figure 1!  

\medskip

$$
\boxit{.5em}{\multiply\hsize2\divide\hsize14 \noindent Yang-Mills\nl
        \noindent Connections} 
\enspace \raise3ex\hbox{$\Leftarrow$}\enspace
\boxit{.5em}{\multiply\hsize2\divide\hsize24 \noindent Moduli\nl
        \noindent Spaces} 
\enspace \raise3ex\hbox{$\Leftarrow$}\enspace
\boxit{.5em}{\multiply\hsize2\divide\hsize16 \noindent Donaldson\nl
        \noindent Invariants} 
\enspace \raise3ex\hbox{$\Leftarrow$}\enspace
\boxit{.5em}{\multiply\hsize2\divide\hsize11 \noindent Seiberg-Witten \nl
        \noindent Theory}
$$

\medskip

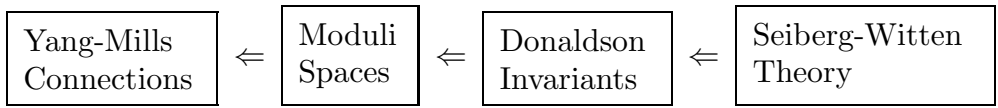
\centerline{Figure 2.  The New Paradigm}

\medskip

We can use SW theory to gain insight into the structure of Donaldson
invariants.  We can use Donaldson invariants to tell us about moduli
spaces.  Finally, we can use moduli spaces to tell us about solutions to
the YM equations.  Perhaps this is not a ``politically correct''
program; nothing points to a classification of 4-manifolds!  But the
structure of moduli spaces can be extremely interesting, so let's
get to work using our new-found tools.

\bigskip

\nd \undertext{\bf Our Results}

\medskip

This paper is an exercise along the lines of Figure 2.  I will take as
given that a manifold has simple type (as all smooth orientable
4-manifolds with $b_+>1$ are believed to have), and see what that says
about the structure of its moduli spaces.  The results are quite
surprising!

Simple type says that the $k+1$st moduli space $\calm_{k+1}$,
intersected with certain varieties, has the homology of 64 copies of
the $k$-th moduli space $\calm_k$.  I will show that the portion of (a
small perturbation of) $\calm_{k+1}$ near the boundary, cut down,
looks like exactly $6$ copies of $\calm_k$, regardless of the topology
or geometry of the underlying 4-manifold.  Simple type thus implies
that the {\it interior} of $\calm_{k+1}$, intersected with certain
varieties, has the homology of 58 copies of $\calm_k$.  This is quite
surprising, since the only known embeddings of $\calm_k$ into
$\calm_{k+1}$ involve Taubes patching, and have images near the
boundary of $\calm_{k+1}$.  Nobody has the slightest idea of what
$\calm_k$ has to do with the interior of $\calm_{k+1}$, yet for all
known 4-manifolds with $b_+>1$, they appear to be closely related.
This is a mystery that warrants further investigation.

Here is an outline for the rest of the paper.  First I review the
definitions of the Donaldson invariants, and of simple type, and state
the result precisely.  Then I sketch the proof, which has three
essential ingredients.  First there is a choice of the geometric
representative of the point class.  Then there is an approximate
formula for the curvature of a connection in the boundary region of
$\calm_{k+1}$.  Finally there is a very naive calculation that
illustrates clearly why the boundary region of $\calm_{k+1}$
contributes exactly 6 copies of $\calm_k$, rather than any other
number.  It takes quite a bit of analysis to thoroughly justify the
approximate formula and the naive calculation; all that can be found
in [GS].  Here I will concentrate on presenting the essential ideas as
clearly as possible, making simplifying assumptions, as needed, along
the way.

\bigskip

\nd \undertext{\bf What Is Simple Type, Anyway?}

\medskip

Let $X$ be an oriented 4-manifold, let $G=SU(2)$ or $SO(3)$ and let
$\calb^*_k$ be the space of irreducible connections (up to gauge
equivalence) on $P_k$, the principal $G$ bundle of instanton number
$k$ over $X$.  Let $\calm_k \subset \calb_k$ be the space of
irreducible connections on $P_k$ with anti-self-dual curvature, modulo
gauge transformations.

Donaldson [D1, D2] defined a map $\mu: H_i(X,Q) \to
H^{4-i}(\calb^*_k,Q)$, $i=$0, 1, 2, 3, whose image freely generates
the rational cohomology of $\calb^*_k$.  Donaldson invariants are
defined by pairing the fundamental class of $\calm_k$ with products of
$\mu$ of the homology classes of $X$, where $k$ is chosen so that the
dimensions match.  Formally, for elements
$[\Sigma_1],\ldots,[\Sigma_n] \in H_*(X)$, we write
$$D([\Sigma_1] \ldots [\Sigma_n]) = \mu([\Sigma_1]) \smile \cdots
\smile \mu([\Sigma_n]) [\calm_k]. \eqno(1) $$ 
Now let $x$ be the point class in $H_0(X)$, and let $\omega$ be any
formal product of classes in $H_*(X)$.  The simple type condition is
that, for all $\omega$,
$$ D(x^2\omega) = 4 D(\omega). \eqno(2) $$

Of course, the ``fundamental class of $\calm_k$'' is usually not well
defined, as $\calm_k$ is typically not compact.  The usual way to make
sense of (1) and (2) is with geometric representatives.  One finds
finite-codimension varieties $V_\Sigma$ in $\calb_k^*$ that are
Poincare dual to $\mu([\Sigma])$.  One then counts points, with sign,
in $V_{\Sigma_1} \cap \cdots \cap V_{\Sigma_n} \cap \calm_k$.  To make
a topological invariant one must show that the number of intersection
points is independent of auxiliary data, such as the metric and the
choice of representatives.  This requires careful analysis of the
bubbling-off phenomena that make $\calm_k$ noncompact.

Unfortunately, $\mu(x)$ is not an integral class in
$H^4(\calb^*)$.  However, $-4\mu(x)$ {\it is} an integral class.  Let
$\nu_1$ and $\nu_2$ be two (generic) geometric representatives of
$-4 \mu(x)$.
The simple type condition can be rewritten as
$$ \#(\calm_{k+1} \cap \nu_1 \cap \nu_2
\cap V_\omega)  = 64 \#(\calm_k \cap V_\omega), \eqno(3) $$
where $\omega$ is an arbitrary
formal product of homology cycles of $X$, and $V_\omega$ is a
geometric representative of $\mu(\omega)$.
Still more formally, one can write
$$ [\calm_{k+1} \cap \nu_1 \cap \nu_2] = 64 [\calm_k].  \eqno(4) $$ I
will show you that, with the right choice of $\nu_i$, the boundary
region of (a perturbation of) $\calm_{k+1}$, intersected with $\nu_1$
and $\nu_2$, looks like 6, not 64, copies of $\calm_k$.

\vfill\eject

\nd \undertext{\bf The Geometric Representative $\nu_p$.}

\medskip

Our first step is to find appropriate geometric representatives
of $\mu(x)$.  The geometric representatives we will use are closely tied
to the notion of reducibility.  An $SU(2)$ connection is said to be
reducible if the gauge group reduces to a proper subgroup of $SU(2)$.
The only proper, nontrivial connected Lie subgroup of $SU(2)$ is
$U(1)$, so the curvature of a reducible connection lives in the Lie
Algebra of $U(1)$, which is 1-dimensional.  Thus all the components of
the curvature must be colinear.  In fact, one can show that, on a
contractible set, an $SU(2)$ connection is reducible if and only if, at
each point, the components of the curvature are colinear.  We
therefore define ``reducible at a point'' to mean that the components of
the curvature are colinear at that point.

Let $p$ be any point in $X$, and let
$$ \nu_p = \{ [A] \in \calb^*_k | F_A^- \hbox{ is reducible at }p \}. 
\eqno(5) $$  
Here $F_A^- = (F_A - *F_A)/2$ is the anti-self-dual part of the
curvature $F_A$.  In [S] I proved the following theorem, which applies
equally well to $SU(2)$ and $SO(3)$ bundles.

\smallskip

\nd {\bf Theorem 1:} {\it $\nu_p$ is a geometric representative of 
$-4 \mu(x)$.}  

\smallskip

\nd {\it Sketch of Proof:} We first need a geometric representative
for the universal $SO(3)$ bundle, and then pull it back to a gauge
theory setting.  Let $V$ be any real
vector space, and let $S_V$ be the Stiefel manifold of linearly
independent triples of vectors in $V$.  $SO(3)$ acts freely on $S_V$,
with a quotient we denote $G_V$.  Topologically, $G_V$ is $\real^6
\times $ the Grassmannian of oriented 3-planes in $V$. If $V$ is
infinite-dimensional, then $S_V$ is contractible, and $S_V \to G_V$ is
the universal $SO(3)$ bundle.  Let $\pi: V \to \real^3$ be any linear
surjection.  If $\dim(V)>7$, then the first
Pontryagin class of the bundle $S_V \to G_V$ is represented by
$$ \nu_\pi = \{ [v_1,v_2,v_3] \in G_V | \pi(v_1), \pi(v_2), \hbox{ and }
\pi(v_3) \hbox{ are colinear.}\} \eqno(6) $$
A proof, using Schubert cycles, may be found in [S], but this result
was almost certainly known to Pontryagin.  

We are now able to construct $\mu$ of the point class.  Let $p$ be a
point on the manifold $X$, let $D$ be a geodesic ball around $p$, let
$\cala_D$ be the $SU(2)$ (or $SO(3)$) connections on $D$ within the
Sobolev space $L^q_k$ (the choice of $q$ and $k$ is not important),
let $\calg^0$ be the gauge transformations in $L^q_{k+1}$ that leave
the fiber at $p$ fixed, and let $\calg$ be all gauge transformations
in $L^q_{k+1}$.  Define $\mu_D(p)$ to be $-{1 \over 4}p_1$ of the
$SO(3)$ bundle $\cala_D/\calg^0 \to \cala_D/\calg$.  $\mu_D(p)$ is a
cohomology class in $H^*(\calb(D))$.  Let $R: \calb(X)
\to \calb(D)$ be the map obtained by restricting connections on a
bundle over $X$ to a bundle over $D$, and define $\mu(x) = R^*
\mu_D(p)$. $\mu(x)$ is then a class in $H^*(\calb(X))$, which turns
out not to depend on the choice of point $p$ or neighborhood $D$.   
For more about the topology of the $\mu$ map, see Chapter 5 of [DK].

Note that, when the gauge group is $SU(2)$, the bundle
$\cala_D/\calg^0 \to \cala_D/\calg$ is a principal $SO(3)$ bundle, not
a principal $SU(2)$ bundle. The reason is that $SU(2)$ does not act
freely on $\cala_D/\calg^0$.  A gauge transformation by $\pm 1$ leaves
a connection fixed, so the typical fiber of our bundle is $SU(2)/Z_2
\sim SO(3)$.  

Although the restriction of the original $SU(2)$ (or $SO(3)$)
bundle to $D$ is trivial, the bundle $\cala_D/\calg^0 \to
\cala_D/\calg$ is highly nontrivial.  Indeed, it is essentially 
the universal $SO(3)$ bundle. The space $\cala_D/\calg^0$ is
isomorphic to the set of connections in radial gauge with respect to
the point $p$.  In such a gauge the connection form $A$ vanishes in
the radial direction but is otherwise unconstrained.  In particular,
$A(p)=0$, and the curvature at $p$, $F_A(p)=dA(p)+ A(p) \wedge A(p) =
dA(p)$, is a linear function of $A$.

Let $V$ be the space of (scalar valued) 1-forms with no radial
component.  A connection in radial gauge is defined by a triple of
elements of $V$, one for each direction in the Lie Algebra.  Deleting
the infinite-codimension set for which these elements are linearly
dependent we get $S_V$. Thus $\mu_D(p)$ is $-1/4 p_1$ of $S_V \to
G_V$, which we have already computed.  Now let $\pi(\alpha)$ be
$d^-(\alpha)$ evaluated at $p$.  $\pi$ is a linear map of $V$ onto the
3-dimensional space of ASD 2-forms at $p$.  $-4\mu_D(p)$ is then
represented by $\nu_\pi$, which is the set of connections on $D$ for
which $F^-(p)$ is reducible.  Pulling $-4 \mu_D(p)$ back by the
restriction map we get the connections on $X$ for which $F_A^-(p)$ is
reducible, i.e. $\nu_p$. \qed

%\vfill\eject

\bigskip

\nd \undertext{\bf The Main Result.}

\medskip

Pick geodesic normal coordinates about some point in $X$, and let $p$
and $q$ be the points $(\pm L, 0,0,0)$, with $L$ small.  Our problem
is to count the points on the left hand side of (3) that lie near the
boundary of $\calm_{k+1}$. This is tantamount to answering the basic
question: {\it Given a connection $A_0 \in \calm_k$, how many ways are
there to glue a small charge-1 bubble onto $A_0$ so as to make the
resulting curvatures reducible at both $p$ and $q$?} This is a local
calculation, and yields an extremely simple answer:

\nd {\bf Theorem 2:} {\it For a generic background connection $A_0$,
and for any $\a \in (0,2)$, there are exactly six ways to glue in a
bubble of size $O(L^\a)$ so as to make $F^-(p)$ and $F^-(q)$ both
reducible.  All six solutions have bubbles of size $O(L^2)$, and all
six have positive orientation.  This answer is independent
of the global topology and geometry of $X$, and in particular is
independent of whether $X$ has simple type.}

\bigskip

\nd \undertext{\bf The Approximate Curvature Formula}

\medskip

Suppose we have a background connection $A_0$, in a radial gauge with
respect to the origin, and glue in a bubble
with center at the origin, size $\l$, and gluing angle $m$ to
get a new connection $A$.  What is the curvature $F_A$ of $A$?
Remarkably, there is an extremely simple approximate formula:
$$ F_A(x) \approx F_{A_0}(x) + F_\std(x), \eqno(7)$$ 
as long as $\l \ll |x| \ll 1$. Here $F_\std$ is the curvature of a
standard $k=1$ instanton, centered at the origin with size $\l$, in a
gauge that is radial, singular at the origin, and regular at $\infty$.
This gauge is not unique; the choice of this gauge is essentially our
gluing angle $m$.  If a bubble is to be glued in at a point $y \ne 0$,
then formula (7) still applies, except that the relevant gauges are
radial with respect to $y$, not to the origin.

The reason for the formula is this.  Let $A_\std$ be the connection
form for the standard instanton.  In the appropriate gauge,
$|A_\std(x)|$ is of order $\l^2/|x|^3$, while $|A_0(x)|= O(|x|)$.  In
the relevant region, the connection form for $A$ is essentially
$A_\std + A_0$, and so the curvature is $F_0 + F_\std + A_\std \wedge
A_0 + A_0 \wedge A_\std$.  The last two terms have norms of order
$\l^2/|x|^2$, and so may be ignored for $|x| \gg \l$.

For the remainder of this paper, we will pretend that (7) is an
equality, rather than an approximation.  The error terms really do
{\it not} matter, although it takes a fair bit of work [GS] to prove it.

Thanks to formula (7), our problem reduces to finding gluing data
such that $F_0 + F_\std$ is reducible at $p$ and $q$.  This involves
two steps: finding what values of $F_\std(p)$ and $F_\std(q)$ are
required, and counting the sets of gluing data that yield those
values.  Both steps require the following notational tool: 

\bigskip

\nd \undertext{\bf Expressing Curvatures as  $3 \times 3$ Real Matrices}

\medskip

Relative to the standard oriented basis of $\Lambda^2_-T^*\real^4$
(namely $\omega_1=dx^0 dx^1 - dx^2 dx^3$, $\omega_2= dx^0 dx^2 - dx^3
dx^1$, $\omega_3=dx^0dx^3 - dx^1 dx^2$), an anti-self-dual curvature
form $F$ has, at each point, 3 Lie-algebra-valued components.  $F$ can
thus be viewed as a triple of 3-vectors, or equivalently a $3 \times 3$
real matrix that we (momentarily) denote $Mat(F)$.  Reducibility at
a point means that this matrix has a rank of 1 (or 0) there.  More
precisely, $Mat(F)$ is constructed as follows.  The first, second and
third columns of $Mat(F)$ are half the $\omega_1$, $\omega_2$ and
$\omega_3$ components of $F$.  The first, second and third entries of
each column refer to the $i$, $j$ and $k$ directions in $su(2)$, the
Lie Algebra of $SU(2)$.  Of course, this construction is dependent on
gauge and a choice of basis for $TX$. A gauge transformation is a
change of basis in $su(2)$, and thus changes $Mat(F)$ by
left-multiplication by an element of $SO(3)$.  A change of basis in
$TX$ changes $Mat(F)$ by right-multiplication by an element of
$SO(3)$.  Thus the singular values of $Mat(F)$, and in particular the
rank of $Mat(F)$, are gauge- and basis-independent.

Now let's compute the matrix of a standard $k=1$ instanton.  Think of
$SU(2)$ as the unit quaternions, with $su(2)$ as the imaginary
quaternions.  The connection form of a standard instanton of scale
size 1, centered at the origin, is $A_{\std_0} = Im(\bar x dx /
(1+|x|^2))$.  The curvature of this connection is
$$F_{\std_0} = {d \bar x dx \over (1 + |x|^2)^2} = 
{2i \omega_1  +  
2j \omega_2  +  
2k \omega_3
\over(1 + |x|^2)^2} \eqno(8) $$
Note that the matrix
$Mat(F_{\std_0})$ is $1/(1 + |x|^2)^2$ times the identity matrix.

Unfortunately, that is in the wrong gauge, in which $A \sim \phi^{-1}
d \phi$ as $|x| \to \infty$, where $\phi(x)=x/|x|$.  We do a gauge
transformation by $\phi^{-1}$, to get a radial gauge in which $A=
O(|x|^{-3})$ as $|x| \to \infty$ (and $A$ is singular at the origin).
We then do a further gauge transformation by an arbitrary constant
$g_0$ to get the most general radial gauge with this property.
$F_\std$ is the curvature form in this gauge.  Since $F_\std =
g_0^{-1} \phi F_{\std_0} \phi^{-1} g_0$, $Mat(F_\std) = \rho(g_0^{-1})
\rho(\phi) Mat(F_{\std_0})$, where $\rho$ is the standard double
covering map from $SU(2)$ to $SO(3)$; the three columns of
$\rho(\phi)$ are $\phi i
\phi^{-1}$,  $\phi j \phi^{-1}$, and $\phi k \phi^{-1}$. The matrix
$\rho(g_0)$ is our gluing angle $m$.

Now suppose that we have a $k=1$ instanton, centered at a point $y$,
with scale size $\lambda$.  The curvature matrix, expressed in the
exterior radial gauge of gluing angle $m$, is then
$$
Mat(F_\std) = {\lambda^2 \over (\lambda^2 + |x-y|^2)^2}\;  
m^{-1} \rho \left ( {x-y \over |x-y|} \right ) \eqno(9) $$
Note that $Mat(F_\std)$ is a positive multiple of an $SO(3)$
matrix.  The multiple is determined by $\l$ and $|x-y|$, while the
$SO(3)$ matrix is determined by $m$ and $(x-y)/|x-y|$.  
  
From now on, we will identify curvatures with their matrices,
and will omit the explicit function ``$Mat$''.

\bigskip

\nd \undertext{\bf A Linear Algebra Lemma}

\medskip

Recall that the singular values $\sigma_1 \ge \sigma_2 \ge \sigma_3
\ge 0$ of a $3 \times 3$ real matrix $M$ are the square roots of the
eigenvalues of $M^T M$.  For $M$ generic, these are distinct and
positive.  The non-generic cases are as follows: Matrices in a
codimension-1 set have $\sigma_3=0$. Matrices in a codimension-2 set
either have $\sigma_1=\sigma_2$ or $\sigma_2=\sigma_3$.  Matrices in a
codimension-4 set have $\sigma_2=\sigma_3=0$; these
matrices have rank 1 or 0. Matries in a codimension-5 set have
$\sigma_1=\sigma_2=\sigma_3$; these are all scalar multiples of
$SO(3)$ matrices.  Only the zero matrix (codimension-9) has
$\sigma_1=\sigma_2=\sigma_3=0$.

\nd {\bf Lemma:} 
{\it Let $P$ be a 3 by 3 real matrix with singular
values $\sigma_1 \ge \sigma_2 \ge \sigma_3 \ge 0$.  If these singular
values are all distinct, then there are exactly two pairs $(s,M)\in
(0,\infty) \times SO(3)$ for which
$P+sM$ has rank 1 (and no pairs $(s,M)$ for which $P+sM=0$).  
In both cases $s =\sigma_2(P)$.  
If exactly two of the singular values of $P$ are the same and
nonzero, then the two solutions $(s,M)$ coalesce to a double root.}

The proof, although straightforward, is not especially enlightening,
so we'll skip it here.  You may want to attempt it as an exercise, or
even assign it to an advanced linear algebra class. One of the many
possible proofs can be found in [GS].

This means that, for any generic background connection $A_0$ and a
generic point $z$, there is one size $s_z$
and two $SO(3)$ matrices $M_1(z)$ and $M_2(z)$ such that $F_0(z)+s_z
M_i(z)$ is reducible.  Note that these are completely determined by
$F_0(z)$, which is a smooth function of $z$.  Since $p$ and $q$ are
separated by a distance $O(L)$, $s_p$ is within $O(L)$ of $s_q$, and
$M_i(p)$ is within $O(L)$ of $M_i(q)$.

\bigskip

\nd \undertext{\bf Solving For The Magnitude Of $F_\std$}

\medskip

According to formula (7), in order to have $F(p)$ and $F(q)$
reducible, we must solve $F_\std(p) = s_p M_i(p)$ and $F_\std(q) = s_q
M_j(q)$.  We do this in two steps, first solving for the magnitudes
and then for the $SO(3)$ matrices.   We will pretend that $s_p=s_q$.
The actual $O(L)$ difference between them complicates the algebra, but
does not make any qualitative difference.  

The condition
for the standard curvature $F_\std$ to have magnitude $s_p$ at $p$ is
$$ {\l^2 \over (|y-p|^2 + \l^2)^2} = s_p, \eqno(10) $$
or equivalently 
$$  \l^2 + |y-p|^2 = \l/\sqrt{s_p}. \eqno(11) $$
Similarly, we need
$$ \l^2 + |y-q|^2 = \l/\sqrt{s_q}. \eqno(12) $$
Since $s_p = s_q$, this implies that $|y-p|=|y-q|$, so $y=(y_0,y_I)$
must lie in the plane half-way between $p$ and $q$.  That is, $y_0=0$.
As long as this is the case, any solution to (11) is also a solution
to (12).

As long as $|y-p| < 1/2\sqrt{s_p}$ there are two solutions to
(11), while for $|y-p| > 1/2\sqrt{s_p}$ there are none.  When $|y-p|<
1/2\sqrt{s_p}$, one solution
has $\l > 1/2 \sqrt{s_p}$, which is not $O(L^\a)$ for $L$
small.  The other solution may qualify as small if $|y-p|$ is small
enough, and, for $|y-p|\ll 1/\sqrt{s_p}$, is approximately $\l =
|y-p|^2 \sqrt{s_p} = (L^2 + |y_I|^2) \sqrt{s_p}$. 

\bigskip

\nd \undertext{\bf Solving For The Gluing Angle}

\medskip

We still have to get the $SO(3)$ matrices right.  By eq.~(9), this means
simultaneously solving the equations 
$$m^{-1} \rho((y\! - \!p)/|y\! -\!p|)= M_i(p) \eqno(13)$$
and 
$$m^{-1} \rho((y\! - \!q)/|y\! - \!q|) = M_j(q) \eqno(14) $$ 
for $m$.
If a solution exists it is obviously unique.  A solution exists if and
only if $\rho((y\! - \!p)/|y\! - \!p|)^{-1}\rho((y\! - \!q)/|y\! -
\!q|) = M_i(p)^{-1} M_j(q)$.  Let 
$$g(y) = {(\bar y -  \bar p)(y  - q)\over |(y  - p)(y - q)|}. \eqno(15)$$
We must count the points on our 3-disk
(of small solutions to (11) and (12)) for which the $SO(3)$-valued
function $\rho(g(y))$ equals $M_i(p)^{-1} M_j(q)$.  Note that
$$g(y) = -I+ 2y_I/L + O((|y_I|/L)^2) \quad \hbox{for }|y_I| \ll L, 
\eqno(16)$$
while 
$$g(y) = I - 2Ly_I/|y_I|^2 + O((L/|y_I|)^2) \quad \hbox{for }|y_I| \gg L. 
\eqno(17)$$ 

Pick a constant $K>0$ and let $R_{K,\alpha}$ be such that
$|y_I|<R_{K,\alpha}$ implies $\l \le K L^\alpha$.  For $L$ small we
have $R_{K,\a}^2 \sim KL^\a/\sqrt{s_m} - L^2 \sim KL^\a/\sqrt{s_m}$,
since $\a < 2$.  $L/R_{K,\alpha}$ is $O(L^{1-\alpha/2})$ and hence
goes to zero as $L \to 0$.  On the 3-disk of admissible $y_I$, the map
$g$ covers all of $SU(2)$ except for a ball of radius $\sim 
L^{1-\alpha/2}$ around the origin.  Since
$\rho$ is a 2-1 map, $\rho(g(y))$ hits all of $SO(3)$ twice, except
for a ball of radius $\sim L^{1-\alpha/2}$ around the
origin, which is only hit once.  The number of solutions to our
problem depends on whether, for small $L$, $M_i^{-1}(p) M_j(q)$ is in
this ball or not.

As $L \to 0$, $M_1(p)^{-1} M_2(q)$ and
$M_2(p)^{-1} M_1(q)$ are bounded away from the identity, but
$M_1(p)^{-1} M_1(q)$ and $M_2(p)^{-1} M_2(q)$ are within $O(L)$ (and
hence within $o(L^{1-\alpha/2})$) of the identity.
Thus we have two sets of parameters $(y,\l,m)$ that give $F_\std(p)=s_p
M_1(p)$ and $F_\std(q)= s_q M_2(q)$, two that give $F_\std(p)= s_p
M_2(p)$ and $F_\std(q)= s_q M_1(q)$, one that gives $F_\std(p) = s_p
M_1(p)$ and $F_\std(q)= s_q M_1(q)$ and one that gives $F_\std(p)= s_p
M_2(p)$ and $F_\std(q) = s_q M_2(q)$.  A total of six solutions in all.

%\bigskip
\vfill\eject

\nd \undertext{\bf Peeking Under The Carpet}

\medskip

If you accept all the simplifying assumptions I have made, then the
proof of Theorem 2 is finished.  I reduced the central
question to counting the solutions to some explicit (and simple!)
algebraic equations, and I not only counted the solutions, but
actually showed you how to find them.  

The skeptical among you, however, may be worried that my simplifying
assumptions are unrealistic, or hide some deep problems.  To ease your
fears, here is a list of the shortcuts I have taken, and how these
issues are dealt with in [GS]. 

I assumed that $s_p=s_q$ when, in fact, $s_p$ and $s_q$ differ by
$O(L)$.  The set of $y$ for which both (11) and (12) can be solved for
$\l$ is typically not a disk in the plane $y_0=0$; rather, it is an
ellipsoid that passes between $p$ and $q$.  The portion of the
ellipsoid that results in $\l$ being small is a topological 3-disk
located between $p$ and $q$.  As $L \to 0$, this topological 3-disk
approaches a geometric 3-disk in the plane $y_0=0$ sufficiently
rapidly that the previous discussion goes through essentially
unchanged.

I showed that there are 6 solutions, but did not show that they all
give intersection number +1.  This involves two steps.  First we show
that the intersection numbers are, by continuity, independent of
$M_i(p)$ and $M_j(q)$.  A similar argument shows that the intersection
numbers for the two solutions for a given $M_i(p)$ and $M_j(q)$ are
equal. We then compute the intersection number at
$(y=0, \l, m=I)$ for $M_i(p)=M_j(q)=I$, which is indeed +1.
 
In order to apply formula (7), $A_0$ must be in radial gauge with
respect to $y$.  However, we computed $M_i(p)$ and $M_j(q)$ from $A_0$
in radial gauge with respect to the origin, not with respect to $y$.
In reality, $M_i(p)$ and $M_j(q)$ should really be viewed as functions of
$y$. The extent of this $y$-dependence can be estimated, and we show
that the derivative of $M_i(p)^{-1}M_j(q)$ with respect to $y$ is too
small to change the count.

Formula (7) is itself an approximation, not an equality.  The error
terms may be treated as a perturbation to $F_0$.  By estimating the
dependence of $(y,\l,m)$ on $F_0$, and the dependence of the error
terms of $(y,\l,m)$, we show that any solution to $F_0(p)+F_\std(p)=$
reducible and $F_0(q) + F_\std(q)=$ reducible can be perturbed to a
solution to $F(p)=$ reducible and $F(q)=$ reducible, and vice versa.
These estimates are, technically, the most difficult part of the whole
problem.

Finally, formula (7) applies not to a true ASD connection, but to a
connection obtained by an explicit grafting formula. The set of such
connections is an $L^2$-small perturbation $\tcalm$ of the true ASD
moduli space $\calm_{k+1}$.  Theorem 2 does
not directly relate $\calm_k$ to $\calm_{k+1}\cap\nu_p\cap\nu_q$.
Rather, it relates $\calm_k$ to $\tcalm\cap\nu_p\cap\nu_q$.  

Ideally, one would like to interpolate from
$\tcalm\cap\nu_p\cap\nu_q$ to $\calm_{k+1}\cap\nu_p\cap\nu_q$.
This is quite difficult, as $\nu_p$ and $\nu_q$ are defined by
pointwise conditions.  I know of no pointwise estimates relating the
curvature of an almost-ASD connection to that of a nearby ASD
connection.  In order to make use of the integral estimates available
in the literature one would have to replace $\nu_p$ and $\nu_q$ by
geometric representatives defined by integral conditions.  While this
is possible (Cliff Taubes once showed me such an extended
representative), it is well beyond the scope of this work.

%\bigskip
\vfill\eject

\nd \undertext{\bf A Differential Forms Approach}

\medskip

There is a quite different approach to measuring the importance of the
boundary region of $\calm_{k+1}$ to simple type.  One can use
differential form representatives of $\mu(\cdot)$, and integrate these
forms over $\calm$ to obtain Donaldson invariants.  In that setting,
our problem is to integrate $\mu_{dR}(p) \wedge \mu_{dR}(q)$ over the
8-dimensional space of gluing parameters.  Here $\mu_{dR}(p)$ is a de
Rham representative of $\mu(x)$ based on connections near $p$, much
the way that $\nu_p$ is a geometric representative.  Let
$B_{\lambda_0}$ be the set of gluing data $(y,\l,m)$ for which $\l <
\l_0$. In [GS] we prove that
$$ \lim_{\l_0 \to 0} \lim_{L \to 0} \int_{B_{\l_0}} \mu_{dR}(p) \wedge
\mu_{dR}(q) = 1/2. \eqno(18)
$$

Surprisingly, this is a different answer than obtained from the
geometric representative calculation (1/8 of what is required for
simple type, as opposed to 6/64).  Moreover, the bulk of the integral
(18) is from $\l$ being of order $L$, while the geometric
representative calculation had all the intersection points having $\l$
of order $L^2$.  This is not a contradiction.  Although the Donaldson
invariants are topological, hence independent of a choice of
representatives, the contribution of the boundary region is geometric,
and can definitely depend on a choice of representatives.  

\bigskip

\nd \undertext{\bf Conclusions}

\medskip

While the two approaches disagree on the exact 
contribution of the boundary region, and on just how close to the
boundary we should consider, they agree on the central theme of this
paper.  Simple type is {\it not} a boundary phenomenon.  Simple type
implies that the features of each moduli space
$\calm_k$ are duplicated  in the structure of the interior of
$\calm_{k+1}$.  This duplication is not at all explained by
our present understanding of moduli spaces; perhaps the explanation
lies in quantum duality.   

%\vfill\eject

\bigskip

\nd \undertext{\bf Acknowledgements}

\medskip

I gratefully acknowledge the collaboration of David
Groisser, and the many helpful discussion I have had with Dan Freed,
Rob Kusner, Tom Parker, Cliff Taubes and Karen Uhlenbeck.  This work
was partially supported by an NSF Mathematical Sciences 
Postdoctoral Fellowship and Texas ARP Grant 003658-037.

% \bigskip

\vfill\eject
\vs .1
\centerline{\undertext{\bf References}}
\vs.2
 
\nd [D1]\ S.K.~Donaldson, 
Connections, cohomology and the intersection forms of four manifolds,
{\it J. Diff. Geom.} 24 (1986), 275--341.
\vs.1 

\nd [D2] \ S.K.~Donaldson,
Polynomial invariants for smooth 4-manifolds,
{\it Topology} 29 (1990), 257--315.
\vs.1

\nd [DK] \ S.K.~Donaldson and P.B.~Kronheimer, 
``The Geometry of Four-Manifolds'', Oxford University Press,
Oxford, 1990.
\vs.1 

\nd[GS] \ D.\ Groisser and L.\ Sadun,
Simple type and the boundary of moduli space,
preprint, 1997.
\vs.1

\nd [S] \ L.\ Sadun, A simple geometric representative of $\mu$ of a
point, {\it Commun. Math. Phys.} 178 (1996), 107--113.
\vs.1 

\nd [W] \ E.\ Witten,
Monopoles and four-manifolds,
{\it Math. Res. Lett.} 1 (1994), 769--796. \vs.1

\bye